\documentclass[a4paper,11pt,english]{article}

\newcommand{\text}{\mbox}
\long\def\TeXButton#1#2{#2}
\newcommand{\proof}{\paragraph{Proof. }}
\def\endproof{\mbox{\ $\Box$}}
\def\limfunc#1{\mathop{\rm #1}}

\usepackage{babel}
\usepackage{oldlfont} 
\newtheorem{theorem}{Theorem}

\newtheorem{proposition}[theorem]{Proposition}

\newtheorem{definition}[theorem]{Definition}
\newtheorem{corollary}[theorem]{Corollary}

\begin{document}

\author{Yuri G. Kondratiev$^{1,2}$ \\Ludwig Streit$^{1,3}$\\Werner Westerkamp$^1$%
\bigskip\  \\ 
$^1$BiBoS, Universit\"at Bielefeld, D 33615 Bielefeld, Germany\\
$^2$Institute of Mathematics, Kiev, Ukraine\\
$^3$CCM, Universidade da Madeira, P 9000 Funchal, Portugal}
\title{A Note on Positive Distributions\\in Gaussian Analysis
\thanks{Published in: Ukrainian Mathematical Journal {\bf 47} No. 5 1995}
}
\date{August 1994\\[10mm]{\large {\it Dedicated to Yu.\ M. Berezansky on the
occasion of his 70th birthday}}}
\maketitle

\begin{abstract}
We describe positive generalized functionals in Gaussian Analysis. We focus
on the large distribution space $({\cal N})^{-1}$. It is shown that a
positive distribution is represented by a measure with specific growth of
its moments. Equivalently this may be replaced by an integrability condition.
\end{abstract}

\section{Introduction}

In recent years Gaussian Analysis and in particular White Noise Analysis
have developed to a useful tool in applied mathematics and mathematical
physics. For a detailed exposition of the theory and for many examples of
applications we refer the reader to the recent monograph \cite{HKPS}.

One of the basic technical ideas in the development of the theory is the use
of dual pairs of spaces of test and generalized functionals. Since the
usefulness of a particular test functional space depends on the application
one has in mind various dual pairs appear in the literature. Here we only
mention two examples which appear very flexible in applications, the
Meyer-Watanabe space ${\cal D}$ and the space of Hida test functionals $(%
{\cal S}).$ Recently new examples of such a dual pair appeared \cite{KoS92}, 
\cite{KLS}. Our interest in the space $({\cal N})^1$ and its dual $({\cal N}%
)^{-1}$ is motivated by the following arguments

\begin{itemize}
\item  It seems to be the natural (and in some sense the only possible){\tt %
\ }construction in generalizations to a non-gaussian setting \cite{AKS}, 
\cite{ADKS}.

\item  The distribution space $({\cal N})^{-1}$ is sufficiently large and
thus contains very singular objects.

\item  An important feature is that it allows the inverse operations for
Wick product and Wick exponential. In fact it allows a powerful 'Wick
calculus', see \cite{KLS}.

\item  This has important applications in the theory of stochastic partial
differential equations \cite{Be93}, \cite{HLOUZ}, \cite{Ok93}, \cite{OkMa}.
\end{itemize}

Motivated by results in finite dimensional distribution theory a natural
question is whether the cone of positive generalized functionals is related
to measures. For the space $({\cal S})^{*}$ of Hida distributions this is
true. Results in this direction can be found in \cite[b]{Ko80a}, \cite{Po87}%
, \cite{BeKo}, \cite{yok90}, \cite{Yok93}, \cite{Lee}.

In this note we prove a similar statement for positive distributions in $(%
{\cal N})^{-1}$. The measures arising in this context are characterized by a
faster growth of their moments compared to the measures related to Hida
distributions. This growth condition can be translated into an equivalent
integrability condition.

The central step in the proof is the construction of the measure by given
positive generalized functional. Since the exponential functionals fail to
be test functionals we are not able to recover the measure from its
characteristic functional using Minlos' theorem. Instead of doing this we
refer to a theorem of Berezansky and Shifrin \cite{BS71} which allows to
construct a unique measure by given moments.

A slight modification of the proof gives a completely analogous result for
the positive distributions in $({\cal N})^{-\beta }$ which will be
formulated without proof.

\section{Gaussian Analysis}

\subsection{Preliminaries}

We start by considering a standard Gelfand triple 
$$
{\cal N}\subset {\cal H\subset N^{\prime }}. 
$$
Here ${\cal H}$ is a real separable Hilbert space with inner product $%
\left\langle \cdot ,\cdot \right\rangle $ and norm $\left| \cdot \right| $
and ${\cal N}$ is a separable nuclear space densely topologically embedded
in ${\cal H}$. For example ${\cal N}$ can be chosen as the Schwartz test
function space ${\cal S}({\bf R})$ for ${\cal H}=L^2({\bf R})$. This
particular choice is the usual one in White Noise Analysis see e.g.\ \cite
{HKPS}. The nuclear space ${\cal N}$ can be represented as projective limit
of a family of Hilbert spaces $\left\{ {\cal H}_p\text{, }p\in {\bf N}%
\right\} $, such that for all $p_1,p_2\in {\bf N}$ there exists $p\in {\bf N}
$ such that ${\cal H}_p\subset {\cal H}_{p_1}$ and ${\cal H}_p\subset {\cal H%
}_{p_2}$ and the embeddings are of Hilbert-Schmidt class (see e.g.\ \cite
{pietsch}). Without loss of generality we suppose that $\forall p\in {\bf N}$%
, $\forall \varphi \in {\cal N}$: $|\varphi |\leq |\varphi |_p$. The dual
space ${\cal N^{\prime }}$ is the inductive limit of the corresponding dual
spaces ${\cal H}_{-p}$. We denote also by $\left\langle \cdot ,\cdot
\right\rangle $ the dual pairings between ${\cal H}_p$ and ${\cal H}_{-p}$
and between ${\cal N}$ and ${\cal N^{\prime }}$ given by the extension of
the inner product on ${\cal H}$. Furthermore $\left| \cdot \right| _{\pm p}$
denote the norms on ${\cal H}_p$ and ${\cal H}_{-p}$ respectively and we
preserve this notation for the norms on the complexifications ${\cal H}_{p,%
{\bf c}}$ and ${\cal H}_{-p,{\bf c}}$ and on tensor powers of these spaces.

The symmetric Fock space over ${\cal H}$ is built up from the symmetric
tensor powers of the complexification ${\cal H}_{{\bf c}}$ : $\Gamma ({\cal %
H)}=\bigoplus\limits_{n=0}^\infty {\cal H}_{{\bf c}}^{\widehat{\otimes }n}$
with inner product given by $\sum_{n=0}^\infty n!\left\langle \overline{%
f^{(n)}},g^{(n)}\right\rangle $. We consider ${\cal C}_\sigma ({\cal %
N^{\prime }})$ the $\sigma $-algebra generated by cylinder sets on ${\cal N}%
^{\prime }$. The canonical Gaussian measure $\mu $ on ${\cal C}_\sigma (%
{\cal N^{\prime }})$ is given by its characteristic functional%
$$
\int_{{\cal N^{\prime }}}\exp (i\left\langle x,\eta \right\rangle )\ {\rm d}%
\mu (x)=\exp (-\frac 12\left| \eta \right| ^2)\text{ , }\eta \in {\cal N} 
$$
via Minlos' theorem (see e.g.\ \cite{BeKo}, \cite{Hi80}). Let us consider
the space of square-integrable complex valued functions with respect to this
measure $L^2(\mu )=L^2({\cal N^{\prime }},{\cal C}_\sigma ({\cal N^{\prime }}%
),\mu )$. The well-known Segal isomorphism between $L^2(\mu )$ and $\Gamma (%
{\cal H)}$ establishes the Wiener-It\^o chaos decomposition of an element $%
\varphi \in L^2(\mu )$ (see e.g.\ \cite{segal}, \cite{pfi2})%
$$
\varphi (x)=\sum_{n=0}^\infty \left\langle :x^{\otimes n}:,\varphi
^{(n)}\right\rangle \text{ .} 
$$
Here $\varphi ^{(n)}\in {\cal H}_{{\bf c}}^{\widehat{\otimes }n}$, $n\in 
{\bf N}$ and $\varphi ^{(0)}\in {\bf C}$ and we have the relation%
$$
\left\| \varphi \right\| _{L^2}^2=\sum_{n=0}^\infty n!\left| \varphi
^{(n)}\right| ^2\quad . 
$$
For the definition of $:x^{\otimes n}:$ we refer the reader to the books of 
\cite{pfi2}, \cite{HKPS}.

\subsection{Test functions and Distributions}

Consider the space ${\cal P}({\cal N^{\prime }})$ of continuous polynomials
on ${\cal N^{\prime }}$, i.e.\ any $\varphi \in {\cal P}({\cal N^{\prime }})$
has the form $\varphi $$(x)=\sum_{n=0}^N\left\langle x^{\otimes n},\tilde
\varphi ^{(n)}\right\rangle $ , $x\in {\cal N}^{\prime }$ , $N\in {\bf N}$
for kernels $\tilde \varphi ^{(n)}\in {\cal N}^{\hat \otimes n}$. It is
well-known that any $\varphi \in {\cal P}({\cal N^{\prime }})$ can be
written as a Wick polynomial i.e.\ $\varphi (x)=\sum_{n=0}^N\left\langle
:x^{\otimes n}:,\varphi ^{(n)}\right\rangle $, $\varphi ^{(n)}\in {\cal N}%
^{\hat \otimes n}$, $N\in {\bf N}$ (see e.g.\ \cite[Ch 2]{BeKo}). To
construct test functions we define for $p,q\in {\bf N}$, $\beta \in \left[
0,1\right] $ the following Hilbertian norm on ${\cal P}({\cal N^{\prime }})$%
\begin{equation}
\label{norms}\left\| \varphi \right\| _{p,q,\beta }^2=\sum_{n=0}^\infty
(n!)^{(1+\beta )}2^{nq}\left| \varphi ^{(n)}\right| _p^2\quad ,\quad \varphi
\in {\cal P}({\cal N^{\prime }})\text{ .} 
\end{equation}
(By $\left| \varphi ^{(0)}\right| _p$ we simply mean the complex modulus for
all $p$.) Then we define $({\cal H}_p)_q^\beta $ to be the completion of $%
{\cal P}({\cal N^{\prime }})$ with respect to $\left\| \cdot \right\|
_{p,q,\beta }$. Or equivalently 
$$
({\cal H}_p)_q^\beta =\left\{ \varphi \in L^2(\mu )\biggm|\left\| \varphi
\right\| _{p,q,\beta }<\infty \right\} \text{ .} 
$$
Finally, the space of test functions $\left( {\cal N}\right) ^\beta $ is
defined to be the projective limit of the spaces $({\cal H}_p)_q^\beta $:%
$$
\left( {\cal N}\right) ^\beta =\bigcap_{p,q\geq 0}\ ({\cal H}_p)_q^\beta 
\text{ . } 
$$
For $0\leq \beta <1$ the corresponding spaces have been studied in \cite
{KoS92} and in the special case of Gaussian product measures all the spaces
for $0\leq \beta \leq 1$ were introduced in \cite{Ko78}. For $\beta =0$ and $%
{\cal N}={\cal S}({\bf R})$ the well-known space $({\cal S})=\left( {\cal S}%
\right) ^0$ of Hida test functions is obtained (e.g.\ \cite{KoSa78}, \cite[b]
{Ko80a}, \cite{KT}, \cite{HKPS}, \cite{BeKo}, \cite{KLPSW}), while in this
work we concentrate on the smallest space $\left( {\cal N}\right) ^1$.

Let $({\cal H}_{-p})_{-q}^{-1}$ be the dual with respect to $L^2(\mu )$ of $(%
{\cal H}_p)_q^1$ and let $\left( {\cal N}\right) ^{-1}$ be the dual with
respect to $L^2(\mu )$ of $\left( {\cal N}\right) ^1$. We denote by $%
\left\langle \!\left\langle \ .\ ,\ .\ \right\rangle \!\right\rangle $ the
corresponding bilinear dual pairing which is given by the extension of the
scalar product on $L^2(\mu )$. We know from general duality theory that%
$$
\left( {\cal N}\right) ^{-1}=\bigcup_{p,q}\ ({\cal H}_{-p})_{-q}^{-1}\quad . 
$$
In particular, we know that every distribution is of finite order i.e.\ for
any $\Phi \in \left( {\cal N}\right) ^{-1}$ there exist $p,q\in {\bf N}$
such that $\Phi \in ({\cal H}_{-p})_{-q}^{-1}$. The chaos decomposition
introduces the following natural decomposition of $\Phi \in \left( {\cal N}%
\right) ^{-1}$. Let $\Phi ^{(n)}\in ({\cal N}_{{\bf c}}^{\prime })^{\widehat{%
\otimes }n}$ be given. Then there is a distribution $\left\langle \Phi
^{(n)},:x^{\otimes n}:\right\rangle $ in $\left( {\cal N}\right) ^{-1}$
acting on $\varphi \in \left( {\cal N}\right) ^1$ as%
$$
\left\langle \!\!\left\langle \left\langle \Phi ^{(n)},:x^{\otimes
n}:\right\rangle ,\varphi \right\rangle \!\!\right\rangle =n!\left\langle
\Phi ^{(n)},\varphi ^{(n)}\right\rangle . 
$$
Any $\Phi \in \left( {\cal N}\right) ^{-1}$ then has a unique decomposition 
$$
\Phi =\sum_{n=0}^\infty \left\langle \Phi ^{(n)},:x^{\otimes
n}:\right\rangle \quad , 
$$
where the sum converges in $\left( {\cal N}\right) ^{-1}$ and we have 
$$
\left\langle \!\left\langle \Phi ,\varphi \right\rangle \!\right\rangle
=\sum\limits_{n=0}^\infty n!\left\langle \Phi ^{(n)},\varphi
^{(n)}\right\rangle \quad ,\ \varphi \in \left( {\cal N}\right) ^1\quad . 
$$
From the definition it is not hard to see that $({\cal H}_{-p})_{-q}^{-1}$
is a Hilbert space with norm%
$$
\left\| \Phi \right\| _{-p,-q,-1}^2=\sum\limits_{n=0}^\infty 2^{-nq}\left|
\Phi ^{(n)}\right| _{-p}^2\quad . 
$$

\noindent {\bf Remark:} Considering also the above mentioned spaces $\left( 
{\cal N}\right) ^\beta $ and their duals $\left( {\cal N}\right) ^{-\beta }$
we have the following chain of spaces 
$$
\left( {\cal N}\right) ^1\subset ...\subset \left( {\cal N}\right) ^\beta
\subset ...\subset \left( {\cal N}\right) =\left( {\cal N}\right) ^0\subset
L^2(\mu )\subset \left( {\cal N}\right) ^{*}\subset ...\subset \left( {\cal N%
}\right) ^{-\beta }\subset ...\subset \left( {\cal N}\right) ^{-1}\quad . 
$$

\subsection{Description of test functions by infinite dimensional holomorphy}

We shall collect some facts from the theory of holomorphic functions in
locally convex topological vector spaces, see e.g.\ \cite{Di}. Let us
introduce the space ${\cal E}({\cal N}_{{\bf c}}^{\prime })$ of complex
valued entire functions on ${\cal N}_{{\bf c}}^{\prime }$. A function $u:%
{\cal N}_{{\bf c}}^{\prime }\rightarrow {\bf C}$ is entire if and only if
for all $y\in {\cal N}_{{\bf c}}^{\prime }$ there exists a sequence of
homogeneous polynomials $\frac 1{n!}{\rm d}^nu(y)(z)$ on ${\cal N}_{{\bf c}%
}^{\prime }$ such that the Taylor decomposition 
$$
u(z+y)=\sum_{n=0}^\infty \frac 1{n!}{\rm d}^nu(y)(z)\text{ ,\quad }z\in 
{\cal N}_{{\bf c}}^{\prime } 
$$
converges uniformly on any neighborhood 
$$
U_{-p,-q}=\left\{ z\in {\cal H}_{-p,{\bf c}}\biggm|\left| z\right|
_{-p}^2<2^q\right\} \quad ,\quad p,q\in {\bf N\quad .} 
$$
We refer to \cite{Di} for the following convenient equivalent definition

\begin{proposition}
A functional $u$ is entire if and only if $u$ is locally bounded and for any 
$y,z\in {\cal N}_{{\bf c}}^{\prime }$ the mapping $\lambda \mapsto
u(y+\lambda z)$, $\lambda \in {\bf C}$ is entire (in the usual sense).
\end{proposition}

We will use a subset ${\cal E}_{\min }^1({\cal N}_{{\bf c}}^{\prime
})\subset {\cal E}({\cal N}_{{\bf c}}^{\prime })$ which consists of all
entire functions of first order of growth and minimal type. I.e.\ $u\in 
{\cal E}({\cal N}_{{\bf c}}^{\prime })$ belongs to ${\cal E}_{\min }^1({\cal %
N}_{{\bf c}}^{\prime })$ iff 
\begin{equation}
\label{mintype}\forall p\in {\bf N\quad }\forall \varepsilon >0\quad \exists
C:\ \left| u(z)\right| \leq Ce^{\varepsilon \left| z\right| _{-p}},\quad
z\in {\cal H}_{-p,{\bf c}}. 
\end{equation}

Now we have introduced the notation to state a theorem proven in \cite{KLS}
which shows that functions from $\left( {\cal N}\right) ^1$ have a pointwise
meaning on ${\cal N}^{\prime }$ and are even (real) analytic on this space.

\begin{theorem}
\label{Description}Any test function in $\left( {\cal N}\right) ^1$ has a
pointwise defined version which has an analytic continuation onto the space $%
{\cal N}_{{\bf c}}^{\prime }$ as an element of ${\cal E}_{\min }^1({\cal N}_{%
{\bf c}}^{\prime })$. Vice versa the restriction of any function in ${\cal E}%
_{\min }^1({\cal N}_{{\bf c}}^{\prime })$ to ${\cal N}^{\prime }$ is in $%
\left( {\cal N}\right) ^1$.
\end{theorem}

\noindent In the rest of the paper we identify any $\varphi \in \left( {\cal %
N}\right) ^1$ with its version in ${\cal E}_{\min }^1({\cal N}_{{\bf c}%
}^{\prime })$. In this sense we may write 
$$
\left( {\cal N}\right) ^1={\cal E}_{\min }^1({\cal N}^{\prime })=\left\{
u\!\mid \!{\cal N}^{\prime }\Bigm|u\in {\cal E}_{\min }^1({\cal N}_{{\bf c}%
}^{\prime })\right\} \text{ .} 
$$
For later use we will prove the following corollary, which gives an explicit
bound of the type (\ref{mintype}) in terms of norms in the spaces $({\cal H}%
_p)_q^1$.

\begin{corollary}
\label{Corexptype} For all $\varphi \in \left( {\cal N}\right) ^1$ and $%
q\geq 0$ we have the following pointwise bound 
\begin{equation}
\label{exptype}\left| \varphi (x)\right| \leq C_{p,\varepsilon \ }\left\|
\varphi \right\| _{p,q,1}e^{\varepsilon \left| x\right| _{-p}}\text{, }x\in 
{\cal H}_{-p}\text{ ,} 
\end{equation}
where $\varepsilon =2^{-\frac q2}$ and%
$$
C_{p,\varepsilon \ }=\int_{{\cal N}^{\prime }}e^{\varepsilon \left| x\right|
_{-p}}\ {\rm d}\mu (x)\text{ .} 
$$
Here $p>0$ is taken such that the embedding $i_0^p:{\cal H}_p{\cal %
\rightarrow H}_0$ is of Hilbert-Schmidt type.
\end{corollary}

\TeXButton{Proof}{\proof}Let us introduce the following function%
$$
w(z)=\sum_{n=0}^\infty (-i)^n\left\langle z^{\otimes n},\varphi
^{(n)}\right\rangle \text{ , }z\in {\cal N}_{{\bf c}}^{\prime } 
$$
using the chaos decomposition $\varphi (x)$ $=\sum_{n=0}^\infty \left\langle
:x^{\otimes n}:,\varphi ^{(n)}\right\rangle $ of $\varphi $. Using the
inequality 
$$
\left| \varphi ^{(n)}\right| _p\leq \frac 1{n!}\varepsilon ^n\left\| \varphi
\right\| _{p,q,1}\text{ , }\varepsilon =2^{-\frac q2} 
$$
we may estimate $\left| w(z)\right| $ for $z\in {\cal H}_{-p}$ as follows%
\begin{eqnarray*}
\left| w(z)\right| &\leq & \sum_{n=0}^\infty \left| \varphi ^{(n)}\right|
_p\left| z\right| _{-p}^n \\
&\leq & \left\| \varphi \right\| _{p,q,1}\sum_{n=0}^\infty \frac
1{n!}\varepsilon ^n\left| z\right| _{-p}^n \\
&=&\left\| \varphi \right\| _{p,q,1}\exp \left( \varepsilon \left| z\right|
_{-p}\right) \text{ .} 
\end{eqnarray*}
To achieve a bound of the type (\ref{exptype}) we use the relation \cite{KLS}%
, \cite{BeKo} 
$$
\varphi (x)=\int_{{\cal N}^{\prime }}w(y+ix)\ {\rm d}\mu (y)\text{ , }x\in 
{\cal N}^{\prime }\text{ .} 
$$
This allows to estimate%
\begin{eqnarray*}
\left| \varphi (x)\right| &\leq &\left\| \varphi \right\| _{p,q,1}\int_{{\cal N}^{\prime }}\exp \left( \varepsilon \left| y+ix\right| _{-p}\right) \ {\rm d}%
\mu (y) \\
&\leq &\left\| \varphi \right\| _{p,q,1}e^{\varepsilon \left| x\right|
_{-p}}\int_{{\cal N}^{\prime }}e^{\varepsilon \left| y\right| _{-p}}\ {\rm d}%
\mu (y)\text{ } 
\end{eqnarray*}
We conclude the proof with the inequality 
$$
C_{p,\varepsilon \ }=\int_{{\cal N}^{\prime }}e^{\varepsilon \left| x\right|
_{-p}}\ {\rm d}\mu (x)\text{ }\leq e^{\frac{\varepsilon ^2}{4\alpha }}\int_{%
{\cal N}^{\prime }}e^{\alpha \left| x\right| _{-p}^2}\ {\rm d}\mu (x)\text{ }
$$
for $\alpha >0$ . If $p>0$ is such that the embedding $i_0^p$ is of
Hilbert-Schmidt type and $\alpha $ is chosen sufficiently small the right
hand integral is finite, see e.g.\ \cite[Fernique's theorem]{Kuo75}. 
\TeXButton{End Proof}{\endproof}

\section{Description of positive distributions}

In this section we will characterize the positive distributions in $({\cal N}%
)^{-1}$. We will prove that the positive distributions can be represented by
measures. They are characterized by a condition of growth of their moments.
In the case of the Hida distribution space $({\cal N})^{*}$ similar
statements can be found in works of Kondratiev \cite[b]{Ko80a} and Yokoi 
\cite{yok90}, \cite{Yok93}, see also \cite{Po87} and \cite{Lee}.

Since test functionals $\varphi \in ({\cal N})^1$ are pointwisely defined
functions in the sense of theorem \ref{Description} it makes sense to define
that $\varphi $ is positive ($\varphi \geq 0$) iff $\varphi (x)\geq 0$ for
all $x\in {\cal N}^{\prime }$. The definition of positivity of test
functionals may be relaxed to $\mu $-a.e. positvity.

\begin{definition}
An element $\Phi \in ({\cal N})^{-1}$ is called positive if for any positive 
$\varphi \in ({\cal N})^{1\text{ }}$we have $\left\langle \!\left\langle
\Phi ,\varphi \right\rangle \!\right\rangle \geq 0$ . The cone of positive
elements in $({\cal N})^{-1}$ is denoted by $({\cal N})_{+}^{-1}$.
\end{definition}

\begin{theorem}
\label{ThPosi} Let $\Phi \in ({\cal N})_{+}^{-1}$ be a positive generalized
function. Then there is a unique (positive) measure $\nu $ on $({\cal N}%
^{\prime }{\cal ,C}_\sigma ({\cal N^{\prime }}){\bf )}$ such that $\forall
\varphi \in ({\cal N}{\bf )}^1$ 
\begin{equation}
\label{posi}\left\langle \!\left\langle \Phi ,\varphi \right\rangle
\!\right\rangle =\int_{{\cal N}^{\prime }}\varphi (x)\ {\rm d}\nu (x) 
\end{equation}
and, moreover, $\exists p\geq 0,\;K,C>0:\forall \xi \in {\cal N},\ n\in {\bf %
N}_0$ 
\begin{equation}
\label{momentgrowth}\left| \int_{{\cal N}^{\prime }}\ \left\langle x,\xi
\right\rangle ^n\ {\rm d}\nu (x)\right| \leq KC^nn!\left| \xi \right| _p^n\ 
\text{ .} 
\end{equation}
Vice versa, any (positive) measure $\nu $ which obeys (\ref{momentgrowth})
defines a positive distribution $\Phi \in ({\cal N})_{+}^{-1}$ by (\ref{posi}%
).
\end{theorem}

\noindent {\bf Remark:} For a given measure $\nu $ the distribution $\Phi $
may be viewed as the generalized Radon-Nikodym derivative $\frac{{\rm d}\nu 
}{{\rm d}\mu }$ of $\nu $ with respect to $\mu $. In fact if $\nu $ is
absolutely continuous with respect to $\mu $ then the usual Radon-Nikodym
derivative coincides with $\Phi .$\medskip\ 

\noindent {\bf Proof of Theorem \ref{ThPosi}. }\\To prove the first part we
define moments of a distribution $\Phi $ and give bounds on their growth.
Using this we construct a measure $\nu $ which is uniquely defined by given
moments\TeXButton{TeX}{\renewcommand{\thefootnote}{\fnsymbol{footnote}}}%
\footnote{%
Since the algebra of exponential functions is not contained in $({\cal N})^1$
(see \cite{KLS}) we can not use Minlos' theorem to construct the measure.
This was the method used in Yokoi's work \cite{yok90}.}. The next step is to
show that any test functional $\varphi \in {\cal (N)}^1$ is integrable with
respect to $\nu $.

Since ${\cal P}\subset {\cal (N)}^1$ we may define moments of a positive
distribution $\Phi \in ({\cal N})^{-1}$ by 
$$
{\rm M}_n(\xi _1,...,\xi _n)=\left\langle \!\!\!\left\langle \Phi ,\
\prod\limits_{j=1}^n\left\langle \cdot ,\xi _j\right\rangle \right\rangle
\!\!\!\right\rangle \ ,\quad \ n\in {\bf N},{\bf \quad }\xi _j\in {\cal N,\ }%
1\leq j\leq n 
$$
$$
{\rm M}_0=\left\langle \!\left\langle \Phi ,\ 1\right\rangle \!\right\rangle
\ \text{.} 
$$
To get estimates on the moments we first assume $\xi _1=...=\xi _n=\xi \in 
{\cal N}.$ Since $\Phi \in ({\cal H}_{-p})_{-q}^{-1}$ for some $p,q>0$ we
may estimate as follows%
$$
\bigg|\left\langle \!\left\langle \Phi ,\left\langle \cdot ,\xi
\right\rangle ^n\right\rangle \!\right\rangle \bigg|\leq \left\| \Phi
\right\| _{-p,-q,-1}\left\| \left\langle \cdot ,\xi \right\rangle ^n\right\|
_{p,q,1}\text{ .} 
$$
To obtain a bound of $\left\| \left\langle \cdot ,\xi \right\rangle
^n\right\| _{p,q,1}$ we use the well known Hermite decomposition%
$$
\left\langle \ \cdot \ ,\xi \right\rangle ^n=\sum_{k=0}^{[\frac n2]}\frac{n!%
}{k!\ (n-2k)!}\left( -\frac 12\left| \xi \right| _0^2\right) ^k:\left\langle
\cdot ,\xi \right\rangle ^{(n-2k)}: 
$$
and the equality%
$$
\left\| :\left\langle \cdot ,\xi \right\rangle ^n:\right\| _{p,q,1}=n!\
2^{\frac 12nq\ }\left| \xi \right| _p^n\text{ .} 
$$
Since Wick powers $:\left\langle \cdot ,\xi \right\rangle ^k:$ , $k\in {\bf N%
}$ are orthogonal in each $({\cal H}_p)_q^1$ we can estimate as follows%
\begin{eqnarray*}
\left\| \left\langle \cdot ,\xi \right\rangle ^n\right\|
_{p,q,1}^2 &=& \sum_{k=0}^{[\frac n2]}\left( \frac{n!}{k!\ (n-2k)!\ 2^k\ }%
\right) ^2\left| \xi \right| _0^{4k}\ \left\| :\left\langle \cdot ,\xi
\right\rangle ^{(n-2k)}:\right\| _{p,q,1}^2 \\
&\leq & (n!)^2\ 2^{nq}\left| \xi \right| _p^{2n}\sum_{k=0}^{[\frac n2]}\frac
1{(k!)^2}2^{-2k(1+q)} \\
&\leq & {\rm I}_0(2^{-q})\ (n!)^2\ 2^{nq}\left| \xi \right| _p^{2n}\text{ ,}
\end{eqnarray*}
where ${\rm I}_0$ is a modified Bessel function of order zero (${\rm I}%
_0(2^{-q})<1.3$ for $q\geq 0$). Via the polarization formula this implies%
$$
\ \left\| \prod\limits_{j=1}^n\left\langle \cdot ,\xi _j\right\rangle
\right\| _{p,q,1}\leq \sqrt{{\rm I}_0(2^{-q})}\left( e\ 2^{\frac q2}\right)
^nn!\prod\limits_{j=1}^n\left| \xi \right| _p\text{ .} 
$$
Then we arrive at 
\begin{equation}
\label{mombound}\Big|{\rm M}_n(\xi _1,...\xi _n)\Big|\leq K\ C^n\
n!\prod\limits_{j=1}^n\left| \xi _j\right| _p 
\end{equation}
with $K=\sqrt{{\rm I}_0(2^{-q})}\cdot \left\| \Phi \right\| _{-p,-q,-1}\quad
,\quad C=e\ 2^{\frac q2}.$

Due to the kernel theorem \cite{BeKo}, \cite{GV} we then have the
representation%
$$
{\rm M}_n(\xi _1,...\xi _n)=\left\langle {\rm M}^{(n)},\xi _1\otimes
...\otimes \xi _n\right\rangle \text{ ,} 
$$
where ${\rm M}^{(n)}\in {\cal N}^{\prime ^{\hat \otimes n}}$. The sequence $%
\left\{ {\rm M}^{(n)},\ n\in {\bf N}_0\right\} $ has the following property
of positivity: for any finite sequence of smooth kernels $\left\{
g^{(n)},n\in {\bf N}\right\} $ (i.e.\ $g^{(n)}\in {\cal N}^{\hat \otimes n}$
and $g^{(n)}=0$\ $\forall \;\,n\geq n_0$ for some $n_0\in {\bf N}${\bf )}
the following inequality is valid 
\begin{equation}
\label{posmon}\sum_{k,j}^{n_0}\left\langle {\rm M}^{(k+j)}\;,g^{(k)}\otimes 
\overline{g^{(j)}}\right\rangle \geq 0\text{ .} 
\end{equation}
This follows from the fact that the left hand side can be written as $%
\left\langle \!\left\langle \Phi ,|\varphi |^2\right\rangle \!\right\rangle $
with%
$$
\varphi (x)=\sum_{n=0}^{n_0}\left\langle x^{\otimes n},g^{(n)}\right\rangle
,\quad x\in {\cal N}^{\prime }\text{ ,} 
$$
which is a smooth polynomial. Following \cite{BS71}, \cite{BeKo}
inequalities (\ref{mombound}) and (\ref{posmon}) are sufficient to ensure
the existence of a uniquely defined measure $\nu $ on $({\cal N}^{\prime },%
{\cal C}_\sigma ({\cal N}^{\prime }))$, such that for any $\varphi \in {\cal %
P}$ we have 
$$
\left\langle \!\left\langle \Phi ,\varphi \right\rangle \!\right\rangle
=\int_{{\cal N}^{\prime }}\varphi (x)\ {\rm d}\nu (x)\text{ .} 
$$

Now we are going to prove the embedding $({\cal N})^1\subset L^1({\cal N}%
^{\prime },{\cal C}_\sigma ({\cal N}^{\prime }),\nu )$. Let $\varphi \in (%
{\cal N})^1$ be given. Since $\varphi $ allows bounds of the type (\ref
{exptype}) we have to find $p^{\prime }$ and $\varepsilon $ such that $\exp
(\varepsilon \left| x\right| _{-p^{\prime }})$ is integrable with respect to 
$\nu $.

Choose $p^{\prime }>p$ such that the embedding $i_p^{p^{\prime }}:{\cal H}%
_{p^{\prime }}\rightarrow {\cal H}_p$ is of Hilbert-Schmidt type.{\bf \ }
Let $\left\{ e_k,\ k\in {\bf N}\right\} \subset {\cal N}$ be an orthonormal
basis in ${\cal H}_{p^{\prime }}$. Then $\left| x\right| _{-p^{\prime
}}^2=\sum\limits_{k=1}^\infty \left\langle x,e_k\right\rangle ^2$, $x\in 
{\cal H}_{-p^{\prime }}$. We will first estimate the moments of even order%
$$
\int_{{\cal N}^{\prime }}\left| x\right| _{-p^{\prime }}^{2n}\ {\rm d}\nu
(x)=\sum\limits_{k_1=1}^\infty \cdots \sum\limits_{k_n=1}^\infty \int_{{\cal %
N}^{\prime }}\left\langle x,e_{k_1}\right\rangle ^2\ \cdots \left\langle
x,e_{k_n}\right\rangle ^2\ {\rm d}\nu (x)\ \text{,} 
$$
where we changed the order of summation and integration by a monotone
convergence argument. Using the bound (\ref{mombound}) we have%
\begin{eqnarray*}
\int_{{\cal N}^{\prime }}\left| x\right| _{-p^{\prime }}^{2n}\ {\rm d}\nu
(x) &\leq & K\ C^{2n}\ (2n)!\sum\limits_{k_1=1}^\infty \cdots
\sum\limits_{k_n=1}^\infty \left| e_{k_1}\right| _p^2\cdots \left|
e_{k_n}\right| _p^2 \\
&=& K\ C^{2n}\ (2n)!\left( \sum\limits_{k=1}^\infty \left| e_k\right|
_p^2\right) ^n \\
&=& K\ \left( C\cdot \left\| i_p^{p^{\prime }}\right\| _{HS}\right) ^{2n}(2n)!
\end{eqnarray*}
because%
$$
\sum\limits_{k=1}^\infty \left| e_k\right| _p^2=\left\| i_p^{p^{\prime
}}\right\| _{HS}^2\text{ .} 
$$
The moments of arbitrary order can now be estimated by the Schwarz inequality%
\begin{eqnarray*}
\int \left| x\right| _{-p^{\prime }}^n\ {\rm d}\nu (x) & \leq & \sqrt{\nu ({\cal N}^{\prime })}\left( \int \left| x\right| _{-p}^{2n}\ {\rm d}\nu (x)\right)
^{\frac 12} \\
&\leq & \sqrt{K\ \cdot \left\| \Phi \right\| _{-p,-q,-1}}\ \left( C\left\|
i_p^{p^{\prime }}\right\| _{HS}\right) ^n\sqrt{(2n)!} \\
&\leq & \sqrt[4]{{\rm I}_0(2^{-q})}\left\| \Phi \right\| _{-p,-q,-1}\ \left( e\
2^{(1+\frac q2)}\left\| i_p^{p^{\prime }}\right\| _{HS}\right) ^nn! 
\end{eqnarray*}
since $(2n)!\leq 4^n(n!)^2$ and $\nu ({\cal N}^{\prime })={\bf E}(\Phi )\leq
\left\| \Phi \right\| _{-p,-q,-1}$. \\Choose $\varepsilon <$ $\left( e\
2^{(1+\frac q2)}\left\| i_p^{p^{\prime }}\right\| _{HS}\right) ^{-1}$ then%
\begin{eqnarray} \label{exponent}
\int e^{\varepsilon \left| x\right| _{-p^{\prime }}}{\rm d}\nu
(x) &=& \sum_{n=0}^\infty \frac{\varepsilon ^n}{n!}\int \left| x\right|
_{-p^{\prime }}^n\ {\rm d}\nu (x) \nonumber \\
&\leq & \sqrt[4]{{\rm I}_0(2^{-q})}\left\| \Phi \right\|
_{-p,-q,-1}\ \sum_{n=0}^\infty \left( \varepsilon \ e\ 2^{(1+\frac
q2)}\left\| i_p^{p^{\prime }}\right\| _{HS}\right) ^n<\infty 
\end{eqnarray}
Hence the first part of the theorem is proven.

Let $\nu $ be a measure on ${\cal N}^{\prime }$ such that (\ref{momentgrowth}%
) holds. Then an argument completely analogous to the first part of the
proof shows that $({\cal N})^1\subset L^1({\cal N}^{\prime },{\cal C}_\sigma
({\cal N}^{\prime }),\nu )$ .

To establish continuity of the linear functional $\varphi \mapsto \int
\varphi (x)\ {\rm d}\nu (x)$, $\varphi \in ({\cal N})^1$ we proceed as
follows. From the first part of the proof we know that (\ref{momentgrowth})
is sufficient to prove that $\exp (\varepsilon \left| x\right| _{-p^{\prime
}})$ is integrable with respect to $\nu $. We only have to choose $p^{\prime
}>p$ such that $i_p^{p^{\prime }}:{\cal H}_{p^{\prime }}\rightarrow {\cal H}%
_p$ is of Hilbert-Schmidt type and to choose $\varepsilon $ sufficiently
small. Because of corollary \ref{Corexptype} there exist $C_{p^{\prime
},\varepsilon \ }>0$ and $q>0$ such that 
$$
\left| \varphi (x)\right| \leq C_{p^{\prime },\varepsilon \ }\left\| \varphi
\right\| _{p,q,1}\exp (\varepsilon \left| x\right| _{-p^{\prime }})\text{ , }%
x\in {\cal H}_{-p^{\prime }}\text{ .} 
$$
Thus 
$$
\left| \int_{{\cal N}^{\prime }}\varphi (x)\ {\rm d}\nu (x)\right| \leq
C_{p^{\prime },\varepsilon \ }\left\| \varphi \right\| _{p,q,1}\int_{{\cal N}%
^{\prime }}\exp (\varepsilon \left| x\right| _{-p^{\prime }})\ {\rm d}\nu
(x) 
$$
Since the integral on the right hand side is finite the measure $\nu $
defines a distribution $\Phi \in ({\cal N})^{-1}$ by%
$$
\left\langle \!\left\langle \Phi ,\varphi \right\rangle \!\right\rangle
=\int_{{\cal N}^{\prime }}\varphi (x)\ {\rm d}\nu (x)\text{ .} 
$$
\TeXButton{End Proof}{\endproof}\bigskip\ 

\noindent {\bf Remark: }From the theorem follows immediately: If $\;\Phi \in
({\cal N})_{+}^{-1}$ then $\Phi $ can be extended to those exponential
functions $\varphi =e^{i\left\langle \cdot ,\xi \right\rangle }$ for which $%
\limfunc{Im}\xi $ is in a neighborhood of zero where $T\Phi \equiv
\left\langle \!\left\langle \Phi ,e^{i\left\langle x,\cdot \right\rangle
}\right\rangle \!\right\rangle =\int_{{\cal N}^{\prime }}e^{i\left\langle
x,\cdot \right\rangle }\ {\rm d}\nu (x)$ is holomorphic. This is due to the
fact that the characteristic function of $\nu $ constructed in theorem \ref
{ThPosi} is analytic in a strip of regularity.\bigskip\ \ 

Using the proof of theorem \ref{ThPosi} we can also show the following
equivalent formulation which uses an integrability condition instead of
moment inequalities:\bigskip\ 

\noindent {\bf Theorem \ref{ThPosi}}$^{\prime }$ \\{\it Let }$\Phi \in (%
{\cal N})_{+}^{-1}$.{\it \ Then there is a unique (positive) measure }$\nu $%
{\it \ on }$({\cal N}^{\prime }{\cal ,C}_\sigma ({\cal N}^{\prime }){\bf )}$%
{\it \ such that }$\forall \varphi \in ({\cal N}{\bf )}^1${\it \ } 
\begin{equation}
\label{posi2}\left\langle \!\left\langle \Phi ,\varphi \right\rangle
\!\right\rangle =\int_{{\cal N}^{\prime }}\varphi (x)\ {\rm d}\nu (x) 
\end{equation}
{\it \ and, moreover, $\exists p^{\prime }>0,\ \exists \varepsilon >0$%
\begin{equation}
\label{adm}\int_{{\cal N}^{\prime }}e^{\varepsilon \left| x\right|
_{-p^{\prime }}}{\rm d}\nu (x)<\infty \ . 
\end{equation}
Vice versa, any (positive) measure }$\nu ${\it \ which obeys (\ref{adm})
defines a positive distribution }$\Phi \in ({\cal N})_{+}^{-1}${\it \ by (%
\ref{posi2}). }

\TeXButton{Proof}{\proof}The first part already is proved in the course of
the proof of theorem \ref{ThPosi}. See equation (\ref{exponent}). The other
direction is trivial since (\ref{adm}) implies a condition of the type (\ref
{momentgrowth}).\TeXButton{End Proof}{\endproof}\bigskip\ 

A straightforward modification in the proof of theorem \ref{ThPosi} gives
the following theorem characterizing $({\cal N})_{+}^{-\beta }$.

\begin{theorem}
Let $\Phi \in ({\cal N})_{+}^{-\beta }$, $\beta \in [0,1)$ be a positive
generalized function. Then there is a unique (positive) measure $\nu $ on $(%
{\cal N}^{\prime }{\cal ,C}_\sigma ({\cal N^{\prime }}){\bf )}$ such that $%
\forall \varphi \in ({\cal N}{\bf )}^\beta $ 
\begin{equation}
\label{posi3}\left\langle \!\left\langle \Phi ,\varphi \right\rangle
\!\right\rangle =\int_{{\cal N}^{\prime }}\varphi (x)\ {\rm d}\nu (x)
\end{equation}
and, moreover, $\exists p\geq 0,\;K,C>0:\forall \xi \in {\cal N},\ n\in {\bf %
N}_0$ 
\begin{equation}
\label{momentgrowth3}\left| \int_{{\cal N}^{\prime }}\ \left\langle x,\xi
\right\rangle ^n\ {\rm d}\nu (x)\right| \leq KC^n(n!)^{\frac{1+\beta }%
2}\left| \xi \right| _p^n\ \text{ .}
\end{equation}
Vice versa, any (positive) measure $\nu $ which obeys (\ref{momentgrowth3})
defines a positive distribution $\Phi \in ({\cal N})_{+}^{-\beta }$ by (\ref
{posi3}).
\end{theorem}

\noindent For $\beta =0$ this result reduces to the well-known theorem in
the case of Hida distributions \cite[b]{Ko80a}, \cite{Yok93}, \cite{Lee}.

This theorem gives another argument (besides the absence of a
characterization theorem) why we do not consider spaces $({\cal N)}^{-\beta
} $\ for $\beta >1$.\ In that case we would have a growth of the moments
with $(n!)^\alpha $\ where $\alpha >1$.$\ $In such a case the uniqueness of
the measure fails to hold (even in a finite dimensional setting). \bigskip\ 

{\bf Acknowledgements: }We thank Prof.\ B. \O ksendal for deepening our
interest in questions connected with $({\cal N})^{-1}$. One of us (W.W.)
gratefully acknowledges financial support of a scholarship from
'Graduiertenf\"orderung des Landes Nordrhein-Westfalen'.

\TeXButton{TeX field}{\footnotesize}


\begin{thebibliography}{HL\O UZ93a}
\bibitem[AKS93]{AKS}  Albeverio, S., Kondratiev, Yu.G. and Streit, L.
(1993), {\it How to generalize White Noise Analysis to Non-Gaussian Spaces}.
In: 'Dynamics of Complex and Irregular Systems'. Eds.: Ph. Blanchard et al.,
World Scientific, 401-130.

\bibitem[ADKS94]{ADKS}  Albeverio, S., Daletzky, Y., Kondratiev, Yu. G.,
Streit, L. (1994), {\it Non-Gaussian infinite dimensional analysis, }%
preprint, 42p.

\bibitem[Be93]{Be93}  Benth, F.E. (1993), {\it A note on population growth
in a crowded stochastic environment, }manuscript, University of Oslo, 21p.

\bibitem[BS71]{BS71}  Berezansky, Yu. M. and Shifrin, S.N. (1971), {\it The
generalized degree symmetric Moment Problem,} Ukrainian Math. J. 23 N3,
247-258.

\bibitem[BeKo]{BeKo}  Berezansky, Yu. M. and Kondratiev, Yu. G. (1988), {\it %
Spectral Methods in Infinite-Dimensional Analysis}, (in Russian), Naukova
Dumka, Kiev. To appear in English in 1994, Kluwer Academic Publishers,
Dordrecht, 800p.

\bibitem[FHSW94]{FHSW}  de Faria, M., Hida, T., Streit, L., and Watanabe,
H., {\it Intersection local times as Generalized White Noise Functionals},
in preparation.

\bibitem[FPS91]{FPS}  de Faria, M., Potthoff, J. and Streit, L. (1991), {\it %
The Feynman integrand as a Hida distribution.} J. Math. Phys. 32, 2123-2127.

\bibitem[Di]{Di}  Dineen, S. (1981), {\it Complex Analysis in Locally Convex
Spaces,} Mathematical Studies 57, North Holland, Amsterdam, 492p.

\bibitem[GV]{GV}  Gel'fand, I.M. and Vilenkin, N.Ya. (1968),{\it \
Generalized Functions}, Vol. IV, Academic Press, New York and London, 384p.

\bibitem[Hi75]{Hi75}  Hida, T. (1975), {\it Analysis of Brownian
Functionals, }Carleton Math. Lecture Notes No. 13, Carleton, 56p.

\bibitem[Hi80]{Hi80}  Hida,T. (1980), {\it Brownian Motion}. Springer, New
York, 325p.

\bibitem[HKPS]{HKPS}  Hida, T., Kuo, H.H., Potthoff, J. and Streit, L.
(1993),{\it \ White Noise. An infinite dimensional calculus}. Kluwer,
Dordrecht, 516p.

\bibitem[HL\O UZ93a]{HL93a}  Holden, H., Lindstr\o m, T., \O ksendal, B.,
Ub\o e, J. and Zhang, T.--S.(1993): {\it Stochastic boundary value problems:
A white noise functional approach}; Probab. Th. Rel. Fields {\bf 95},
391--419.

\bibitem[HL\O UZ93b]{HLOUZ}  Holden, H., Lindstr\o m, T., \O ksendal, B.,
Ub\o e, J. and Zhang, T.--S.(1993): {\it The pressure equation for fluid
flow in a stochastic medium,} Preprint 18p.

\bibitem[Ko78]{Ko78}  Kondratiev, Yu.G. (1978), {\it Generalized functions
in problems of infinite dimensional analysis}. Ph.D. thesis, Kiev
University, 175p.

\bibitem[Ko80a]{Ko80a}  Kondratiev, Yu.G. (1980), {\it Spaces of entire
functions of an infinite number of variables, connected with the rigging of
a Fock space. }In: 'Spectral Analysis of Differential Operators.' Math.
Inst. Acad. Sci. Ukrainian SSR, p. 18-37. English translation: Selecta Math.
Sovietica {\bf 10} (1991), 165-180.

\bibitem[Ko80b]{Ko80b}  Kondratiev, Yu.G. (1980), {\it Nuclear spaces of
entire functions in problems of infinite dimensional analysis.} Soviet Math.
Dokl. {\bf 22}, 588-592.

\bibitem[KLPSW94]{KLPSW}  Kondratiev, Yu.G., Leukert, P., Potthoff, J.,
Streit, L., Westerkamp, W. (1994), {\it Generalized Functionals in Gaussian
Spaces - the Characterization Theorem Revisited. }Manuskripte 175/94, Uni
Mannheim, 15p.

\bibitem[KLS94]{KLS}  Kondratiev, Yu.G., Leukert, P., Streit, L. (1994),{\it %
\ Wick Calculus in Gaussian Analysis,} BiBoS preprint 637, 23p, to appear in
Acta Applicandae Mathematicae.

\bibitem[KoSa76]{KoSa76}  Kondratiev, Yu.G. and Samoilenko, Yu.S. (1976),%
{\it \ Integral representation of generalized positive definite kernels of
an infinite number of variables}. Soviet Math. Dokl. 17, 517-521.

\bibitem[KoSa78]{KoSa78}  Kondratiev, Yu.G. and Samoilenko, Yu.S. (1978), 
{\it Spaces of trial and generalized functions of an infinite number of
variables}, Rep. Math. Phys. {\bf 14}, No.3, 325-350.

\bibitem[KoS93]{KoS92}  Kondratiev, Yu.G. and Streit, L. (1993), {\it Spaces
of White Noise distributions: Constructions, Descriptions, Applications.} I.
Rep. Math. Phys. {\bf 33}, 341-366.

\bibitem[KT80]{KT}  Kubo, I. and Takenaka, S. (1980a), {\it Calculus on
Gaussian white noise I, II.} Proc. Japan Acad. 56, 376-380 and 411-416.

\bibitem[Kuo75]{Kuo75}  Kuo, H.-H.(1975), {\it Gaussian Measures in Banach
Spaces.} LNM 463, Springer, New York,168p.

\bibitem[Kuo92]{Kuo92}  Kuo, H.-H. (1992), {\it Lectures on white noise
analysis}. Soochow J. Math. {\bf 18}, 229-300.

\bibitem[LLSW94]{LLSW}  Lascheck, A., Leukert, P., Streit, L., Westerkamp,
W. (1994), {\it More about Donsker's Delta Function}. Soochow J. Math. {\bf %
20} No.3, 401-418.

\bibitem[Lee91]{Lee}  Lee, Y.-J. (1991), {\it Analytic Version of Test
Functionals, Fourier Transform and a Characterization of Measures in White
Noise Calculus. }J. Funct. Anal. {\bf 100}, 359-380.

\bibitem[Lu]{Lu}  Lukacs, E. (1970), {\it Characteristic Functions}, 2nd
edition, Griffin, London, 350p.

\bibitem[\O k93]{Ok93}  \O ksendal, B. (1993), {\it Some mathematical models
for population growth in a stochastic environment}, preprint, 32p.

\bibitem[\O k94]{OkMa}  \O ksendal, B. (1994),{\it \ Stochastic Partial
Differential Equations and Applications to Hydrodynamics. }In: 'Stochastic
Analysis and Applications in Physics.' Ed.: A.I. Cardoso et al.; Kluwer,
Dordrecht, in print.

\bibitem[Pi]{pietsch}  Pietsch, A (1969), {\it Nukleare Lokal Konvexe Raeume}%
, Berlin, Akademie Verlag, 362p.

\bibitem[Po87]{Po87}  Potthoff, J. (1987), {\it On positive generalized
functionals. }J. Funct. Anal. {\bf 74}, 81-95.

\bibitem[Po91]{Po91}  Potthoff, J. (1991), {\it Introduction to white noise
analysis}. In: 'Control Theory, Stochastic Analysis and Applications.' Eds.:
S. Chen, J. Yong; Singapore, World Scientific, 241-256.

\bibitem[Po92]{Po92}  Potthoff, J. (1992), {\it White noise methods for
stochastic partial differential equations.} In: 'Stochastic Partial
Differential Equations and Their Applications.' Eds.: B.L. Rozovskii, R.B.
Sowers; Berlin, Heidelberg, New York, Springer, 28-42.

\bibitem[Po94]{Po94}  Potthoff, J.: {\it White noise approach to parabolic
stochastic differential equations. }In: 'Stochastic Analysis and
Applications in Physics.' Ed.: A.I. Cardoso et al.; Kluwer, Dordrecht, in
print.

\bibitem[PS91]{PS}  Potthoff, J. and Streit, L. (1991), {\it A
characterization of Hida distributions.} J. Funct. Anal. {\bf 101}, 212-229.

\bibitem[PS93]{PSfi}  Potthoff, J. and Streit, L. (1993), {\it Invariant
states on random and quantum fields: }$\phi -${\it bounds and white noise
analysis}. J. Funct. Anal. {\bf 101}, 295-311.

\bibitem[Si]{pfi2}  Simon, B. (1974), {\it The }$P(\Phi )_2${\it \ Euclidean
(Quantum) Field Theory}. Princeton University Press, Princeton, 392p.

\bibitem[Se56]{segal}  Segal, I. (1956), {\it Tensor algebras over Hilbert
spaces}. Trans. Amer. Math. Soc. {\bf 81}, 106-134.

\bibitem[SW93]{SW}  Streit, L. and Westerkamp, W. (1993),{\it \ A
generalization of the characterization theorem for generalized functionals
of White Noise}. In: 'Dynamics of Complex and Irregular Systems.' Eds.: Ph.
Blanchard et al.; World Scientific, 175-187.

\bibitem[Yok90]{yok90}  Yokoi, Y.(1990), {\it Positive generalized white
noise functionals}. Hiroshima Math. J. {\bf 20}, 137-157.

\bibitem[Yok93]{Yok93}  Yokoi, Y. (1993), {\it Simple setting for white
noise calculus using Bargmann space and Gauss transform.} Preprint, 33p.
\end{thebibliography}
\end{document}